\documentclass[12pt]{amsart}

\usepackage{enumitem}
\usepackage{psfrag}
\usepackage{graphicx}
\usepackage{pinlabel}
\usepackage{graphicx}
\usepackage{amsmath}
\usepackage{amssymb}
\usepackage{amscd}
\usepackage{autobreak}
\usepackage{picinpar}
\usepackage{times}
\usepackage{pb-diagram}
\usepackage{graphicx}
\usepackage{wrapfig}
\usepackage{xspace}
\usepackage{hyperref}
\usepackage{url}
\usepackage{caption}
\usepackage{subcaption}
\usepackage{fancyhdr}
\usepackage{overpic}
\usepackage{color}
\usepackage[square,comma,sort&compress,numbers]{natbib}
\usepackage{latexsym}
\usepackage{mathrsfs}
\usepackage{amsfonts}
\usepackage{hyperref}
\usepackage{amsthm}
\usepackage{appendix}
\usepackage{pdfsync}
\DeclareMathOperator{\area}{Area}
\DeclareMathOperator{\arccoth}{arccoth}
\DeclareMathOperator{\arctanh}{arctanh}

\theoremstyle{definition}
\newtheorem{theorem}{Theorem}[section]
\newtheorem{lemma}[theorem]{Lemma}

\newtheorem{definition}[theorem]{Definition}

\newtheorem{remark}[theorem]{Remark}
\newtheorem*{theorem*}{Theorem}
\def\qed{\hfill{Q.E.D.}\smallskip}

\begin{document}

\title{\bf Combinatorial curvature flows for generalized hyperbolic circle packings}
\author{Te Ba, Chao Zheng}

\date{\today}

\address{School of Mathematics, Hunan University, Changsha, 410082, P.R.China}
 \email{batexu@hnu.edu.cn}

\address{School of Mathematics and Statistics, Wuhan University, Wuhan 430072, P.R. China}
\email{czheng@whu.edu.cn}


%52C26 Circle packings and discrete conformal geometry
%57Q15 Triangulating manifolds
%53E99 (Geometric evolution equations) None of the above, but in this section
%57M50 General geometric structures on low-dimensional manifolds
%51M10 Hyperbolic and elliptic geometries (general) and generalizations


\begin{abstract}
Generalized circle packings were introduced in \cite{Ba-Hu-Sun} as a generalization of tangential circle packings in hyperbolic background geometry.
In this paper, we introduce the combinatorial Calabi flow, fractional combinatorial Calabi flow and combinatorial $p$-th Calabi flow for generalized hyperbolic circle packings.  
We establish several equivalent conditions regarding the longtime behaviors of these flows. 
This provides effective algorithms for finding the generalized circle packings with prescribed total geodesic curvatures.
\end{abstract}

\maketitle

\section{Introduction}
\subsection{Background}
The notion of circle patterns was proposed in the work of Thurston \cite{thurston} as a significant tool to study the hyperbolic structure on 3-manifolds.
He introduced hyperbolic circle patterns on a triangulated surface with prescribed intersection angles.
The induced polyhedral metric may produce conical singularities at the vertices.
The classical discrete Gaussian curvature is introduced to describe the singularities at the vertices, 
which is defined by the difference of $2\pi$ and the cone angle at the vertices.
%Geometric flows are significant tools in the study of geometry and topology.
%The most famous work related to geometric flows is that of Hamilton \cite{ham}, which introduced Ricci flow to solve the Poincare conjecture. 
%Given a triangulated surface, Thurston \cite[Chapter 13]{thurston} introduced the circle packing metric, which is a type of piecewise flat cone metric with singularities at the vertices.
Motivated by the work of Hamilton \cite{ham}, Chow-Luo \cite{chow-luo} introduced the combinatorial Ricci flow on closed triangulated surfaces, which is a discrete analogue of Hamilton's Ricci flow.
Under some combinatorial conditions, they proved that the combinatorial Ricci flow exists for all time and converges exponentially fast to Thurston's circle patterns on closed triangulated surfaces both in Euclidean and hyperbolic background geometry. 
%Inspired by the work of Calabi \cite{calabi1,calabi2} and Chow-Luo \cite{chow-luo}, Ge \cite{ge} intruduced combinatorial Calabi flow
Since then, combinatorial curvature flows became an important approach for finding geometric structures on low-dimensional manifolds. 
See, for instance,  combinatorial Yamabe flow \cite{luo}, combinatorial Calabi flow \cite{Ge1,Ge2} and fractional combinatorial Calabi flow \cite{Wu-Xu}. 
%combinatorial $p$-th Calabi flow \cite{L-Z} 

Recently, a new geometric data called ``total geodesic curvature'' was introduced by Nie \cite{nie} to measure the singularities of circle patterns in spherical background geometry.
Nie \cite{nie} provided the existence and rigidity results for spherical circle patterns with respect to the total geodesic curvature.
Motivated by Nie's work, Ba-Hu-Sun \cite{Ba-Hu-Sun} investigated the existence and rigidity of the generalized hyperbolic circle packing (the intersection angle of two circles is zero) with respect to the total geodesic curvature. 
To search the generalized circle packings with prescribed total geodesic curvature,
they \cite{Ba-Hu-Sun} further introduced the combinatorial Ricci flow and proved the solution of the combinatorial Ricci flow exists for all time and converges exponentially.
In this paper, we introduce the combinatorial Calabi flow, fractional combinatorial Calabi flow and combinatorial $p$-th Calabi flow for generalized hyperbolic circle packings. 
We further prove the longtime existence and convergence for the solutions of these combinatorial curvature flows.

\subsection{Set up}
We begin by introducing pseudo ideal triangulation on surfaces, 
which is introduced in \cite{penner} and generalized in \cite{guo-luo,Ba-Hu-Sun}. 
Let $(S,T)$ be a connected closed surface with a triangulation $T$. 
Let $V$, $E$, $F$ be the vertex, edge and face set of $T$. For simplicity of notations, we use one index to denote a vertex, two indices to denote an edge ($ij$ is the arc on $S$ joining $i$, $j$). For each $i\in V$, we use $U(i)$ to denote a small open regular neighborhood of $i$. 
We define \[N(I):=\cup_{i\in I}U(i)\] for each $I\subset V$. 
Suppose $I_1,I_2\subset V$ satisfying $I_1\cap I_2=\emptyset$. 
Set \[S_{I_1,I_2}=S\setminus(N(I_1)\cup I_2).\]
Then $S_{I_1,I_2}$ is a connected surface with $n\geq 0$ boundary components and $m\geq 0$ punctures, where $\vert I_1\vert=n$, $\vert I_2\vert=m$.
The intersection
\[T_{I_1,I_2}=\{\sigma\cap S_{I_1,I_2}\vert\sigma\in T\}\]
is called the pseudo ideal triangulation of $S_{I_1,I_2}$. We use $(S_{I_1,I_2},T_{I_1,I_2})$ to denote the surface $S_{I_1,I_2}$ with a pseudo ideal triangulation $T_{I_1,I_2}$.
%\[T_{I_1,I_2}:=T\cap S_{I_1,I_2}\]
%is called a pseudo ideal triangulation of $S_{I_1,I_2}$. 
The intersections
%\[\widetilde{E}=E\cap\widetilde{S},\quad\widetilde{F}=F\cap\widetilde{S}.\]
\[E_{I_1,I_2}:=\{ij\cap S_{I_1,I_2}\vert ij\in E\},\quad F_{I_1,I_2}:=\{ijk\cap S_{I_1,I_2}\vert ijk\in F\}\]
are called the edge and face set of $T_{I_1,I_2}$.
The intersection of a face of $F_{I_1,I_2}$ and $\partial S_{I_1,I_2}$ is called a $B$-arc.

A generalized hyperbolic circle packing metric on $(S_{I_1,I_2},T_{I_1,I_2})$ is a map $k:V\to\mathbb{R}_{+}$ satisfying
\begin{itemize}
\item[($a$)] $k(i)<1$ if $i\in I_1$,
\item[($b$)] $k(i)=1$ if $i\in I_2$,
\item[($c$)] $k(i)>1$ if $i\in I_3$,
\end{itemize}
where $I_3=V\setminus(I_1\cup I_2)$.
The geometry of $(S_{I_1,I_2},T_{I_1,I_2})$ is determiend as follows:
\begin{itemize}
\item[($i$)] The length of edges of $E_{I_1,I_2}$ is defined by $d:E_{I_1,I_2}\to\mathbb{R}_{+}$, where
\[d(ij)=\left\{
\begin{aligned}
&\arctanh k(i)+\arctanh k(j),&i,j\in I_1,\\
&\arccoth k(i)+\arccoth k(j),&i,j\in I_3,\\
&\arctanh k(i)+\arccoth k(j),&i\in I_1, j\in I_3,\\
&+\infty,&\ i\ \text{or}\ j\in I_2,
\end{aligned}
\right.\]
\item[($ii$)] Each angle at the endpoints of the $B$-arcs is defined to be $\pi/2$.
\end{itemize}
It is proved in \cite{Ba-Hu-Sun} that the side lengths of $B$-arcs and angles of each face can be uniquely determined by the generalized hyperbolic circle packing metric on $(S_{I_1,I_2},T_{I_1,I_2})$ and  ($i$), ($ii$). 
%See \cite{Ba-Hu-Sun} for the detailed discussion.

Let us provide a brief introduction to the geometric meaning of $d(ij)$. 
If $i,j\in I_1$, $d(ij)$ is the distance between of axis of two hypercycles with curvature $k(i)$, $k(j)$. 
If $i,j\in I_3$, $d(ij)$ is the distance between the centers of two circles with curvature $k(i)$, $k(j)$. 
If $i\in I_1$ and $j\in I_3$, $d(ij)$ is the distance between the center of the circle with curvature $k(i)$ and the axis of the hypercycle with curvature $k(j)$. 
If $i\ \text{or}\ j\in I_2$, $d(ij)$ is the distance between the center of the circle with curvature $k(i)=1$ (a horocycle) to the center or axis of a circle, or a horocycle, or a hypercycle with curvature k(v), which is $+\infty$. 

\begin{figure}[htbp]
\centering
\captionsetup{width=0.90\linewidth}
\includegraphics[scale=0.25]{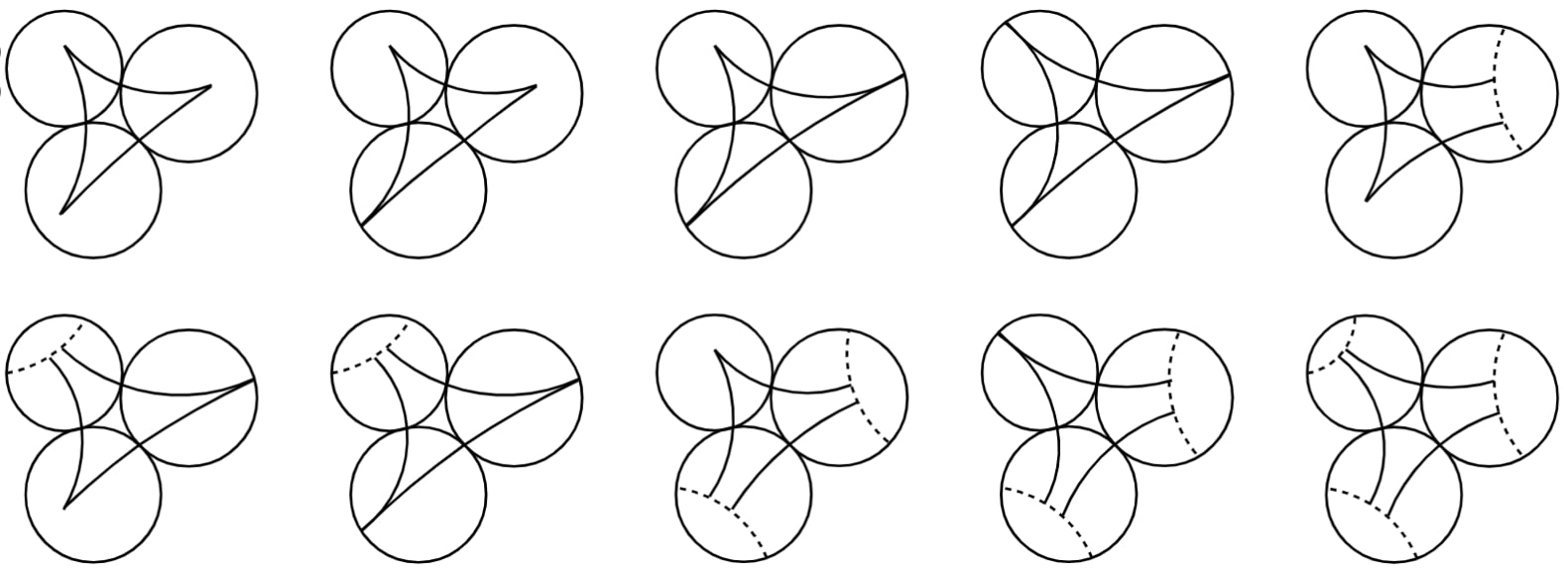}
\caption{\small Three-circle configurations}
\label{fig4}
\end{figure}

Suppose $k:V\to\mathbb{R}_{+}$ is a generalized hyperbolic circle packing metric on $(S_{I_1,I_2},T_{I_1,I_2})$. Each $f\in F_{I_1,I_2}$ can be embedded into three mutually tangent hyperbolic circles (including horocycles and hypercycles).
Here we cite \cite[Figure 4]{Ba-Hu-Sun} as an explanation, as shown in Figure \ref{fig4}.
Then there exists a hyperbolic circle packing (with possibly horocycles or hypercycles) on $S_{I_1,I_2}$ induced by $k$. Let $C_v$ be the circle of this packing which centered at $v$. The total geodesic curvature of $k$ at $v\in V$ is defined as the total geodesic curvature of $C_v$. It can be calculated by
\[L_v=\int_{C_v}k(v)ds=l(v)k(v),\]
where $l(v)$ is the length of $C_v$. Note that each circle is not necessarily to be embedded in $\mathbb{H}^2$ because $(S_{I_1,I_2},T_{I_1,I_2})$ may exists singularities at the vertices or edges when assigned with a generalized hyperbolic circle packing metric.
%The total geodesic curvature was first introduced in the work of Nie \cite{nie} as an important tool to study the existence and rigidity of the circle patterns in spherical background geometry. 
\begin{theorem}[\cite{Ba-Hu-Sun}, Theorem 1.2]\label{Thm: existence}
Let $(S,T)$ be a connected closed surface with the vertex, face set $V$, $F$. 
Let $F_I$ be the set of faces having at least one vertex in $I$ for subset $I\subset V$.
Then there exists $I_1,I_2\subset V$ and a generalized hyperbolic circle packing metric on $(S_{I_1,I_2},T_{I_1,I_2})$ having the total geodesic curvature $L_1,\cdots,L_{\vert V\vert}$ on each vertex if and only if $(L_1,\cdots,L_{\vert V\vert})\in\Omega$, where \begin{equation}\label{Eq: key}
\Omega=\left\{(L_1,\cdots,L_{\vert V\vert})\in\mathbb{R}^{\vert V\vert}_{+}\left\vert\right.\sum\nolimits_{i=1}^{\vert V\vert} L_i<\pi\vert F_I\vert\ \ \text{for each}\ I\subset
 V\right\}.\end{equation}
Moreover, the choice of $I_1,I_2$ and generalized hyperbolic circle packing metric is unique if it exists.
\end{theorem}

%\begin{equation}\label{Eq: CRF}
%\frac{dK_i}{dt}=\widehat{L}_i-L_i
%\end{equation}
%and proved the following theorem.
%
%\begin{theorem}[\cite{Ba-Hu-Sun}, Theorem 1.3]\label{Thm: CRF}
%Let $F_I$ be the set of faces having at least one vertex in $I$ for subset $I\subseteq V$.
%Suppose $\{\widehat{L}_i\}_{i\in V}$ is the prescribed total geodesic curvature.
%The solution of the combinatorial curvature flow (\ref{Eq: CRF}) exists for all time.
%Furthermore, it converges exponentially fast if and only if $\{\widehat{L}_i\}_{i\in V}$ satisfies the condition (\ref{Eq: key}).
%\end{theorem}

\subsection{Main results}
Motivated by Ge's work on combinatorial Calabi flow  \cite{Ge1,Ge2}, Wu-Xu's work on fractional combinatorial Calabi flow \cite{Wu-Xu} and Lin-Zhang's work on combinatorial $p$-th Calabi flow \cite{L-Z}, 
we introduce the following combinatorial Calabi flow, fractional combinatorial Calabi flow and combinatorial $p$-th Calabi flow for generalized hyperbolic circle packing metric on $(S_{I_1,I_2},T_{I_1,I_2})$.
%and prove the convergence of these flows.
Set $K_i=\ln k_i$. 
\begin{definition}
Let $\widehat{L}\in \mathbb{R}_{+}^{|V|}$ be a given function defined on $V$.
The combinatorial Calabi flow for generalized hyperbolic circle packing metrics on $(S_{I_1,I_2},T_{I_1,I_2})$ is defined to be
\begin{eqnarray}\label{Eq: CF}
\begin{cases}
\frac{dK_i}{dt}=\Delta(L-\widehat{L})_i,\\
K_i(0)=K_0,
\end{cases}
\end{eqnarray}
where $\Delta$ is the discrete Laplace operator defined by
\begin{equation}\label{Eq: Laplace}
\Delta f_i
=-\sum_{j=1}^{|V|}\frac{\partial L_i}{\partial K_j}f_j
\end{equation}
for any function $f: V\rightarrow \mathbb{R}$.
\end{definition}

Set
\begin{equation*}%\label{Eq: matrix Lambda}
\Lambda=(\Lambda_{ij})_{|V|\times |V|}
=\frac{\partial(L_1,...,L_{|V|})}{\partial(K_1,...,K_{|V|})}.
\end{equation*}
The equation (\ref{Eq: Laplace}) implies $\Delta_{\mathcal{T}}=-\Lambda$.
By Lemma \ref{Lem: matrix}, the matrix $\Lambda$ is symmetric and positive definite on $\mathbb{R}^{|V|}$.
There exists an orthonormal matrix $Q$ such that
\begin{equation*}
\Lambda=Q^T\cdot \text{diag}\{\lambda_1,...,\lambda_{|V|}\}\cdot Q,
\end{equation*}
where $\lambda_1,...,\lambda_{|V|}$ are non-negative eigenvalues of the matrix $\Lambda$.
For any $s\in \mathbb{R}$, the $2s$-th order fractional discrete Laplace operator $\Delta^s$ is defined to be
\begin{equation}\label{Eq: fractional Laplace}
\Delta^s=-\Lambda^s=-Q^T\cdot \text{diag}\{\lambda^s_1,...,\lambda^s_{|V|}\}\cdot Q.
\end{equation}
Therefore, the fractional discrete Laplace operator $\Delta^s$ is negative definite on $\mathbb{R}^{|V|}$.
Specially, if $s=0$, then $\Delta^s$ is reduced to the minus identity operator;
if $s=1$, then $\Delta^s$ is reduced to the discrete Laplace operator $\Delta=-\Lambda=-(\frac{\partial L_i}{\partial K_j})_{|V|\times |V|}$.

\begin{definition}\label{Def: FCF}
Let $\widehat{L}\in \mathbb{R}_{+}^{|V|}$ be a given function defined on $V$.
The fractional combinatorial Calabi flow for generalized hyperbolic circle packing metrics on $(S_{I_1,I_2},T_{I_1,I_2})$ is defined to be
\begin{eqnarray}\label{Eq: FCF}
\begin{cases}
\frac{dK_i}{dt}=\Delta^s(L-\widehat{L})_i,\\
K_i(0)=K_0,
\end{cases}
\end{eqnarray}
where $\Delta^s$ is the fractional discrete Laplace operator defined by (\ref{Eq: fractional Laplace}).
\end{definition}

\begin{remark}
If $s=0$, the fractional combinatorial Calabi flow (\ref{Eq: FCF}) is reduced to the combinatorial Ricci flow introduced by Ba-Hu-Sun \cite{Ba-Hu-Sun}. 
If $s=1$, the fractional combinatorial Calabi flow (\ref{Eq: FCF}) is reduced to the combinatorial Calabi flow (\ref{Eq: CF}).
\end{remark}

By Lemma \ref{Lem: matrix}, we have
\begin{equation*}
\Delta f_i
=-\Lambda f_i
=\sum_{j\sim i}(-B_{ij})(f_j-f_i)+A_if_i,
\end{equation*}
where
$B_{ij}=\frac{\partial L^{jk}_i}{\partial K_j}+\frac{\partial L^{jl}_i}{\partial K_j}$ defined by (\ref{Eq: B}) and
$A_i=\frac{\partial}{\partial K_i}\left(\sum_{ijk}\text{Area}(\Omega_{ijk})\right)$ defined by (\ref{Eq: A}).
For any $p>1$, we define the discrete $p$-th Laplace operator $\Delta_{p}$ for generalized hyperbolic circle packing metrics by the following formula
\begin{equation}\label{Eq: P-Laplace}
\Delta_pf_i
=\sum_{j\sim i}(-B_{ij})|f_j-f_i|^{p-2}(f_j-f_i),
\end{equation}
where $f: V\rightarrow \mathbb{R}$ is a function.

\begin{definition}\label{Def: PCF}
Let $\widehat{L}\in \mathbb{R}_{+}^{|V|}$ be a given function defined on $V$.
The combinatorial $p$-th Calabi flow for generalized hyperbolic circle packing metrics on $(S_{I_1,I_2},T_{I_1,I_2})$ is defined to be
\begin{eqnarray}\label{Eq: PCF}
\begin{cases}
\frac{dK_i}{dt}=(\Delta_p+A_i)(L-\widehat{L})_i,\\
K_i(0)=K_0,
\end{cases}
\end{eqnarray}
where $\Delta_p$ is the discrete $p$-th Laplace operator defined by (\ref{Eq: P-Laplace}).
\end{definition}

\begin{remark}
If $p=2$, then the discrete $p$-th Laplace operator (\ref{Eq: P-Laplace}) is reduced to the discrete Laplace operator (\ref{Eq: Laplace}) and hence the combinatorial $p$-th Calabi flow (\ref{Eq: PCF}) is reduced to the combinatorial Calabi flow (\ref{Eq: CF}).
\end{remark}

The main result of this paper is as follows, which gives
the longtime existence and convergence for the solutions of the combinatorial Calabi flow (\ref{Eq: CF}), the fractional combinatorial Calabi flow (\ref{Eq: FCF}) and the combinatorial $p$-th Calabi flow (\ref{Eq: PCF}) for generalized hyperbolic circle packing metrics on $(S_{I_1,I_2},T_{I_1,I_2})$.

\begin{theorem}\label{Thm: main 1}
Let $\widehat{L}\in \mathbb{R}_{+}^{|V|}$ be a given function defined on $V$.
The following statements are equivalent:
\begin{description}
\item[(1)] $\{\widehat{L}_i\}_{i\in V}\in\Omega$, where $\Omega$ is defined by (\ref{Eq: key});
\item[(2)] The solution of the combinatorial Calabi flow (\ref{Eq: CF}) exists for all time and converges exponentially fast to a unique generalized circle packing metric with the total geodesic curvature $\widehat{L}$;
\item[(3)] The solution of the fractional combinatorial Calabi flow (\ref{Eq: FCF}) exists for all time and converges exponentially fast to a unique generalized circle packing metric with the total geodesic curvature $\widehat{L}$;
\item[(4)] The solution of the combinatorial $p$-th Calabi flow (\ref{Eq: PCF}) exists for all time and converges to a unique generalized circle packing metric with the total geodesic curvature $\widehat{L}$.
\end{description}
%The following statements are equivalent:
%\begin{description}
%  \item[(i)] $\{\widehat{L}_i\}_{i\in V}$ satisfies the condition (\ref{Eq: key}).
%  \item[(ii)]  The solution of the combinatorial Calabi flow (\ref{Eq: CF}) converges for any initial data;
% \item[(iii)] The solution of the fractional combinatorial Calabi flow (\ref{Eq: FCF}) converges for any initial data;
% \item[(iv)] The solution of the combinatorial p-th Calabi flow (\ref{Eq: PCF}) converges for any initial data;
%\end{description}
%Furthermore, if one of the above statements holds, then the solution of the combinatorial Calabi flow (\ref{Eq: CF}) or the fractional combinatorial Calabi flow (\ref{Eq: FCF}) exists for all time and converges exponentially fast or
%the solution of the combinatorial p-th Calabi flow (\ref{Eq: PCF}) exists for all time and converges fast
%to a unique generalized circle packing metric with the total geodesic curvature $\widehat{L}_i$ at $i\in V$ respectively.
\end{theorem}

\begin{remark}
Different from the combinatorial Calabi flow ($p=2$),
we can not get the exponential convergence for the solution of the combinatorial $p$-th Calabi flow for $p>1,p\neq2$.
\end{remark}

%\subsection{Organization of the paper}
%The paper is organized as follows.
%In Section \ref{Sec: matrix lemma}, we give some useful lemmas.
%In Section \ref{Sec: proof of Theorem}, we prove Theorem \ref{Thm: main 1}.

\
\

\noindent\textbf{Acknowledgements}\\[8pt]
The first author is supported by NSF of China (No.11631010). 
The authors would like to thank Xu Xu and Ze Zhou for helpful discussions.

\section{Some useful Lemmas}\label{Sec: matrix lemma}
Let $C_i, C_j, C_k$ be three mutually tangent hyperbolic circles (with possibly horocycles or hypercycles).  
Denote $\Omega_{ijk}$ as the region enclosed by three arcs between tangency points of $C_i$, $C_j$ and $C_k$.
Denote $l_i$ as the length of the arc between two points of tangency of $C_i$.
Set $L_i=l_ik_i$ and $K_i=\ln k_i$.
 
\begin{lemma}[\cite{Ba-Hu-Sun}]\label{Lem: key}
Let $L_{i}^{jk}$, $K_i$ and $\Omega_{ijk}$ be defined as above. Then
\begin{description}
\vspace{0.2cm}
\item[($i$)]
$\frac{\partial L_{i}^{jk}}{\partial K_j}=\frac{\partial L_{j}^{ik}}{\partial K_i}<0.$
\vspace{0.3cm}
\item[($ii$)]
$\frac{\partial\area(\Omega_{ijk})}{\partial K_i}<0.$
\vspace{0.3cm}
\item[($iii$)] $\frac{\partial L_{i}^{jk}}{\partial K_i}>0.$
\end{description}
\end{lemma}
For any two adjacent faces $ijk$ and $ijl$ sharing a common edge $ij$, we set
\begin{equation}\label{Eq: B}
B_{ij}=
\frac{\partial L^{jk}_i}{\partial K_j}
+\frac{\partial L^{jl}_i}{\partial K_j},
\end{equation}
\begin{equation}\label{Eq: A}
A_i
=\frac{\partial}{\partial K_i}\left(\sum_{ijk}
\text{Area}(\Omega_{ijk})\right).
\end{equation}
By Lemma \ref{Lem: key}, we have
$A_i<0$ and $B_{ij}<0$.
Set $\Lambda_A=\text{diag}\{A_1,...,A_{|V|}\}$ and 
$\Lambda_B=((\Lambda_B)_{ij})_{|V|\times |V|}$, where 
\begin{eqnarray*}
(\Lambda_B)_{ij}=
\begin{cases}
-\sum_{k\sim i}B_{ik},  &{j=i},\\
B_{ij}, &{j\sim i},\\
0, &{j\nsim i, j\neq i}.
\end{cases}
\end{eqnarray*}
Then
$\Lambda_A$ is negative definite.
For any $x\in \mathbb{R}^{|V|}$, we have
\begin{equation*}
x^T\Lambda_Bx=\sum_{i,j=1}^{|V|}(\Lambda_B)_{ij}x_ix_j
=\sum_{i\sim j}(B_{ij}x_ix_j)-\sum_{i\sim j}x_i^2B_{ij}=-\frac{1}{2}\sum_{i\sim j}B_{ij}(x_i-x_j)^2\geq 0,
\end{equation*}
which implies $\Lambda_B$ is positive semi-definite.

\begin{lemma}\label{Lem: matrix}
The matrix $\Lambda=\frac{\partial(L_1,...,L_{|V|})}{\partial(K_1,...,K_{|V|})}$ could be decomposed to be
\begin{equation*}
\Lambda=-\Lambda_A+\Lambda_B.
\end{equation*}
As a result, the matrix $\Lambda$ is symmetric and positive definite on $\mathbb{R}^{|V|}$.
\end{lemma}
\proof

By Lemma \ref{Lem: key},
$\frac{\partial L^{jk}_i}{\partial K_j}
=\frac{\partial L^{ik}_j}{\partial K_i}$ and then $\frac{\partial L_i}{\partial K_j}=\frac{\partial L_j}{\partial K_i}$.

\begin{description}
  \item[(1)] If $j\nsim i$ and $j\neq i$, then $\frac{\partial L_i}{\partial K_j}=0.$
  \item[(2)] If $j\sim i$, then
\begin{equation*}
\frac{\partial L_i}{\partial K_j}
=\frac{\partial(\sum_{ijk}L^{jk}_i)}{\partial K_j}
=\sum_{ijk}\frac{\partial L^{jk}_i}{\partial K_j}
=\frac{\partial L^{jk}_i}{\partial K_j}+\frac{\partial L^{jl}_i}{\partial K_j}.
\end{equation*}
Then
$\frac{\partial L_i}{\partial K_j}=B_{ij}$ by (\ref{Eq: B}).
\item[(3)] If $j=i$, then
\begin{equation*}
\begin{aligned}
\frac{\partial L_i}{\partial K_i}
=&\sum_{ijk}\frac{\partial \left(\pi-L^{ik}_j-L^{ij}_k-\text{Area}(\Omega_{ijk})\right)}{\partial K_i}\\
=&-\sum_{ijk}\left(\frac{\partial L^{jk}_i}{\partial K_j}+\frac{\partial L^{jk}_i}{\partial K_k}\right)
-\sum_{ijk}\frac{\partial \text{Area}(\Omega_{ijk})}{\partial K_i}\\
=&-\sum_{j\sim i}\left(\frac{\partial L^{jk}_i}{\partial K_j}+\frac{\partial L^{jl}_i}{\partial K_j}\right)-\frac{\partial}{\partial K_i}\left(\sum_{ijk}\text{Area}(\Omega_{ijk})\right),
\end{aligned}
\end{equation*}
\end{description}
where the first equality is due to following formula obtained by Ba-Hu-Sun (\cite{Ba-Hu-Sun}, Lemma 2.10) 
\begin{equation*}
\mathrm{Area}(\Omega_{ijk})=\pi-L^{jk}_i-L^{ki}_j-L^{ij}_k.
\end{equation*}
Then 
$\frac{\partial L_i}{\partial K_i}=-A_i-\sum_{j\sim i}B_{ij}$.
Therefore, $\Lambda=-\Lambda_A+\Lambda_B$.
\qed

\section{The proof of Theorem \ref{Thm: main 1}}
\label{Sec: proof of Theorem}
We divide Theorem \ref{Thm: main 1} into three theorems and prove them respectively.

\begin{theorem}\label{Thm: main 2}
Let $\widehat{L}\in \mathbb{R}_{+}^{|V|}$ be a given function defined on $V$.
If the solution of the combinatorial Calabi flow (\ref{Eq: CF}) converges for any initial data, then $\{\widehat{L}_i\}_{i\in V}\in\Omega$.
Furthermore, if $\{\widehat{L}_i\}_{i\in V}\in\Omega$,
then the solution of the combinatorial Calabi flow (\ref{Eq: CF}) exists for all time and converges exponentially fast to a unique generalized circle packing metric with the total geodesic curvature $\widehat{L}$.
\end{theorem}
\proof
As $t\rightarrow +\infty$, the solution $K(t)$ of the combinatorial Calabi flow (\ref{Eq: CF}) converges to $\widehat{K}$.
By the $C^1$-smoothness of $L$, we have 
$L(\widehat{K})=\lim_{t\rightarrow +\infty}L(K(t))$.
By the mean value theorem, there exists a sequence $\xi_n\in(n,n+1)$ such that
\begin{equation*}
K_i(n+1)-K_i(n)=K'_i(\xi_n)=\Delta(L(K(\xi_n))
-\widehat{L})_i\rightarrow 0,\ \text{as}\ n\rightarrow +\infty.
\end{equation*}
Combining with Lemma \ref{Lem: matrix}, we have $L_{i}(\widehat{K})=\lim_{n\rightarrow +\infty}L_{i}(K(\xi_n))=\widehat{L}_i$ for all $i\in V$ .
Then $\widehat{K}$ is a generalized circle packing metric with the total geodesic curvature $\widehat{L}$.
By Theorem \ref{Thm: existence}, $\{\widehat{L}_i\}_{i\in V}\in\Omega$.

Conversely, if $\{\widehat{L}_i\}_{i\in V}\in\Omega$, 
there exists a unique $\widehat{K}$ with the total geodesic curvature $\widehat{L}$ by Theorem \ref{Thm: existence}.
The following function
\begin{equation*}
\mathcal{E}(K)=\int^K_{\widehat{K}}
\sum_{i=1}^{|V|}(L_i-\widehat{L}_i)dK_i.
\end{equation*}
is well-defined and strictly convex on $\mathbb{R}^{|V|}$ by Lemma \ref{Lem: matrix}.
Furthermore, $\mathcal{E}(\widehat{K})=0,\ \nabla\mathcal{E}(\widehat{K})=0$ and $\mathrm{Hess}\ \mathcal{E}>0$. 
This implies $\lim_{||K||\rightarrow +\infty}\mathcal{E}(K)=+\infty$.
Hence, $\mathcal{E}(K)$ is proper and $0=\mathcal{E}(\widehat{K})\leq \mathcal{E}(K)$.
By direct calculations, we have
\begin{equation*}
\frac{d\mathcal{E}(K(t))}{dt}
=\sum^{|V|}_{i=1}\frac{\partial \mathcal{E}}{\partial K_i}\frac{dK_i}{dt}
=\sum^{|V|}_{i=1}(L-\widehat{L})_i\Delta(L-\widehat{L})_i
=-(L-\widehat{L})^T\cdot \Lambda\cdot (L-\widehat{L})\leq0
\end{equation*}
by Lemma \ref{Lem: matrix},
which implies $0\leq\mathcal{E}(K(t))\leq \mathcal{E}(K(0))$.
Thus the solution $\{K(t)\}$ of the combinatorial Calabi flow (\ref{Eq: CF}) lies in a compact subset of $\mathbb{R}^{|V|}$,
which implies the solution of the combinatorial Calabi flow (\ref{Eq: CF}) exists for all time.

By Lemma \ref{Lem: key}, the matrix $\Lambda^2$ is strictly positive definite on $\mathbb{R}^{|V|}$.
By the continuity of the eigenvalues of $\Lambda^2$,
there exists $\lambda_0>0$ such that the non-zero eigenvalues $\lambda$ of $\Lambda^2$ satisfy $\lambda>\lambda_0$ along the combinatorial Calabi flow (\ref{Eq: CF}).
Therefore, for the combinatorial Calabi energy
\begin{equation}\label{Eq: 3}
\mathcal{C}(K)
=||L-\widehat{L}||^2
=\sum_{i=1}^{|V|}(L_i-\widehat{L}_i)^2,
\end{equation}
we have
\begin{equation*}
\frac{d\mathcal{C}(K(t))}{dt}
=\sum_{i=1}^{|V|}\frac{\partial\mathcal{C}}{\partial K_i}\frac{d K_i}{dt}
=-2(L-\widehat{L})^T\cdot\Lambda^2\cdot(L-\widehat{L})
\leq -2\lambda_0\mathcal{C}(K(t)),
\end{equation*}
which implies
$\mathcal{C}(K(t))
\leq e^{-2\lambda_0t}\mathcal{C}(K(0))$.
Then
\begin{equation*}
||K(t)-\widehat{K}||^2
\leq C_1||L(t)-\widehat{L}||^2
\leq C_1e^{-2\lambda_0t}||L(0)-\widehat{L}||^2
\leq C_2e^{-2\lambda_0t}
\end{equation*}
for some positive constants $C_1,C_2$.
\qed

\begin{theorem}\label{Thm: main 3}
Let $\widehat{L}\in \mathbb{R}_{+}^{|V|}$ be a given function defined on $V$.
If the solution of the fractional combinatorial Calabi flow (\ref{Eq: FCF}) converges for any initial data, then $\{\widehat{L}_i\}_{i\in V}\in\Omega$.
Furthermore, if $\{\widehat{L}_i\}_{i\in V}\in\Omega$,
then the solution of the fractional combinatorial Calabi flow (\ref{Eq: FCF}) exists for all time and converges exponentially fast to a unique generalized circle packing metric with the total geodesic curvature $\widehat{L}$.
\end{theorem}
\proof
The first part follows from $\Delta^s$ is a negative definite on $\mathbb{R}^{|V|}$.
If $\{\widehat{L}_i\}_{i\in V}\in\Omega$, there exists a unique $\widehat{K}$ with the total geodesic curvature $\widehat{L}$ by Theorem \ref{Thm: existence}.
By direct calculations, we have
\begin{equation*}
\frac{d\mathcal{E}(K(t))}{dt}
=\sum^{|V|}_{i=1}\frac{\partial \mathcal{E}}{\partial K_i}\frac{dK_i}{dt}
=\sum^{|V|}_{i=1}(L-\widehat{L})_i\Delta^s(L-\widehat{L})_i
=-(L-\widehat{L})^T\cdot \Lambda^s\cdot (L-\widehat{L})\leq0
\end{equation*}
by Lemma \ref{Lem: matrix},
which implies $0\leq\mathcal{E}(K(t))\leq \mathcal{E}(K(0))$.
Combining with the properness of $\mathcal{E}$,
the solution $\{K(t)\}$ of the fractional combinatorial Calabi flow (\ref{Eq: FCF}) lies in a compact subset of $\mathbb{R}^{|V|}$,
which implies the solution of the fractional combinatorial Calabi flow (\ref{Eq: FCF}) exists for all time and $\mathcal{E}(K(t))$ converges.
There exists a sequence $\xi_n\in(n,n+1)$ such that as $n\rightarrow +\infty$,
\begin{equation*}
\begin{aligned}
&\mathcal{E}(K(n+1))-\mathcal{E}(K(n))
=(\mathcal{E}(K(t))'|_{\xi_n}
=\nabla \mathcal{E}\cdot\frac{dK_i}{dt}|_{\xi_n}\\
=&\sum^{|V|}_{i=1}(L-\widehat{L})_i\Delta^s(L-\widehat{L})_i|_{\xi_n}
=-(L-\widehat{L})^T\cdot \Lambda^s\cdot (L-\widehat{L})|_{\xi_n}\rightarrow 0.
\end{aligned}
\end{equation*}
Then $\lim_{n\rightarrow +\infty}L_i(K(\xi_n))=\widehat{L}_i
=L_i(\widehat{K})$ for all $i\in V$.
By $\{K(t)\}\subset\subset \mathbb{R}^{|V|}$, there exists $K^*\in \mathbb{R}^{|V|}$ and a subsequence of $\{K(\xi_n)\}$, still denoted as $\{K(\xi_n)\}$ for simplicity, such that $\lim_{n\rightarrow \infty}K(\xi_n)=K^*$, which implies
$L_i(K^*)=\lim_{n\rightarrow +\infty}L_i(K(\xi_n))
=L_i(\widehat{K})$.
This further implies $K^*=\widehat{K}$ by Theorem \ref{Thm: existence}.
Therefore, $\lim_{n\rightarrow \infty}K(\xi_n)=\widehat{K}$.

Set $\Gamma(u)=\Delta^s (L-\widehat{L})$,
then $D\Gamma|_{K=\widehat{K}}$ is negative definite,
which implies that $\widehat{K}$ is a local attractor of (\ref{Eq: FCF}).
The conclusion follows from Lyapunov Stability Theorem (\cite{Pontryagin}, Chapter 5). 

\qed

\begin{remark}
One can also use the combinatorial Calabi energy (\ref{Eq: 3}) to prove the exponential convergence of the solution $K(t)$ of the fractional combinatorial Calabi flow (\ref{Eq: FCF}), which is similar to the proof of Theorem \ref{Thm: main 2}.
\end{remark}

\begin{theorem}\label{Thm: main 4}
Let $\widehat{L}\in \mathbb{R}_{+}^{|V|}$ be a given function defined on $V$.
If the solution of the combinatorial $p$-th Calabi flow (\ref{Eq: PCF}) converges for any initial data, then $\{\widehat{L}_i\}_{i\in V}\in\Omega$.
Furthermore, if $\{\widehat{L}_i\}_{i\in V}\in\Omega$,
then the solution of the combinatorial $p$-th Calabi flow (\ref{Eq: PCF}) exists for all time and converges to a unique generalized circle packing metric with the total geodesic curvature $\widehat{L}$.
\end{theorem}
\proof
Suppose the solution $K(t)$ of the combinatorial $p$-th Calabi flow (\ref{Eq: PCF}) converges to $\widehat{K}$ as $t\rightarrow +\infty$, then
$L(\widehat{K})=\lim_{t\rightarrow +\infty}L(K(t))$ by the $C^1$-smoothness of $L$.
Furthermore, there exists a sequence $\xi_n\in(n,n+1)$ such that
\begin{equation*}
K_i(n+1)-K_i(n)=K'_i(\xi_n)=(\Delta_p+A_i)(L(K(\xi_n))
-\widehat{L})_i\rightarrow 0,\ \text{as}\ n\rightarrow +\infty.
\end{equation*}
Set $\widetilde{L}
=\lim_{n\rightarrow +\infty}(L(K(\xi_n))-\widehat{L})
=L(\widehat{K})-\widehat{L}$,
then 
\begin{equation}\label{Eq: 6}
\lim_{n\rightarrow +\infty}
\sum_{i=1}^{|V|}\widetilde{L}_i(\Delta_p+A_i)\widetilde{L}_i
=0.
\end{equation}
Since $A_i<0$ by (\ref{Eq: A}),
then
\begin{equation}\label{Eq: 7}
\sum_{i=1}^{|V|}\widetilde{L}_iA_i\widetilde{L}_i\leq0.
\end{equation}
By the following formula obtained by Lin-Zhang (\cite{L-Z}, Lemma 5.5)
\begin{equation*}
\sum_{i=1}^{|V|}f_i\Delta_{p}f_i
=\frac{1}{2}\sum_{i=1}^{|V|}\sum_{j\sim i}B_{ij}|f_j-f_i|^p
\end{equation*}
for any $f: V\rightarrow \mathbb{R}$,
we have
\begin{equation}\label{Eq: 8}
\sum_{i=1}^{|V|}\widetilde{L}_i\Delta_p\widetilde{L}_i\leq0
\end{equation}
by $B_{ij}<0$ in (\ref{Eq: B}).
Combining (\ref{Eq: 6}), (\ref{Eq: 7}) and (\ref{Eq: 8}),
we have $\widetilde{L}=0$, i.e.,
$L_{i}(\widehat{K})=\widehat{L}_i$ for all $i\in V$.
By Theorem \ref{Thm: existence}, $\{\widehat{L}_i\}_{i\in V}\in\Omega$.

Conversely, if $\{\widehat{L}_i\}_{i\in V}\in\Omega$, there exists a unique $\widehat{K}$ with the total geodesic curvature $\widehat{L}$ by Theorem \ref{Thm: existence}.
By direct calculations, we have
\begin{equation*}
\begin{aligned}
\frac{d\mathcal{E}(K(t))}{dt}
=&\sum^{|V|}_{i=1}\frac{\partial \mathcal{E}}{\partial K_i}\frac{dK_i}{dt}
=\sum^{|V|}_{i=1}(L-\widehat{L})_i(\Delta_p+A_i)(L-\widehat{L})_i\\
=&\frac{1}{2}\sum_{i=1}^{|V|}\sum_{j\sim i}B_{ij}|(L-\widehat{L})_i-(L-\widehat{L})_j|^p
+(L-\widehat{L})^T\cdot \Lambda_A\cdot (L-\widehat{L})\\
\leq&0
\end{aligned}
\end{equation*}
by $A_i<0$ and $B_{ij}<0$,
which implies $0\leq\mathcal{E}(K(t))\leq \mathcal{E}(K(0))$.
Combining with the properness of $\mathcal{E}$,
the solution $\{K(t)\}$ of the combinatorial $p$-th Calabi flow (\ref{Eq: PCF}) lies in a compact subset of $\mathbb{R}^{|V|}$,
which implies the solution of the combinatorial $p$-th Calabi flow (\ref{Eq: PCF}) exists for all time and $\mathcal{E}(K(t))$ converges.
There exists a sequence $\xi_n\in(n,n+1)$ such that as $n\rightarrow +\infty$,
\begin{equation*}
\begin{aligned}
\mathcal{E}(K(n+1))-\mathcal{E}(K(n))
=&(\mathcal{E}(K(t))'|_{\xi_n}
=\nabla \mathcal{E}\cdot\frac{dK_i}{dt}|_{\xi_n}\\
=&\sum^{|V|}_{i=1}(L-\widehat{L})_i(\Delta_p+A_i)(L-\widehat{L})_i|_{\xi_n}\\
=&\frac{1}{2}\sum_{i=1}^{|V|}\sum_{j\sim i}B_{ij}|(L-\widehat{L})_i-(L-\widehat{L})_j|^p|_{\xi_n}
+(L-\widehat{L})^T\cdot \Lambda_A\cdot (L-\widehat{L})|_{\xi_n}\\
\rightarrow& 0.
\end{aligned}
\end{equation*}
Then $\lim_{n\rightarrow +\infty}L_i(K(\xi_n))=\widehat{L}_i
=L_i(\widehat{K})$ for all $i\in V$.
By $\{K(t)\}\subset\subset \mathbb{R}^{|V|}$, there exists $K^*\in \mathbb{R}^{|V|}$ and  a convergent subsequence $\{u(\xi_{n_k})\}$ of $\{K(\xi_n)\}$ such that $\lim_{k\rightarrow \infty}K(\xi_{n_k})=K^*$.
Then 
$L_i(K^*)=\lim_{k\rightarrow +\infty}L_i(K(\xi_{n_k}))
=L_i(\widehat{K})$.
This implies $K^*=\widehat{K}$ by Theorem \ref{Thm: existence}.
Therefore, $\lim_{k\rightarrow \infty}K(\xi_{n_k})=\widehat{K}$.

We use Lin-Zhang's trick in \cite{L-Z} to prove  $\lim_{t\rightarrow\infty}K(t)=\widehat{K}$.
Suppose otherwise, there exists $\delta>0$ and $t_n\rightarrow +\infty$ such that
$|K(t_n)-\widehat{K}|>\delta$.
This implies $\{K(\xi_n)\}\subseteq \mathbb{R}^{|V|}\backslash B(\widehat{K},\delta)$,
where $B(\widehat{K},\delta)$ is a ball centered at $\widehat{K}$ with radius $\delta$.
Hence, for any $K\in \mathbb{R}^{|V|}\backslash B(\widehat{K},\delta)$, $\mathcal{E}(u)\geq C>0$.
Then $\mathcal{E}(K(\xi_n))\geq C>0$.
Since $\mathcal{E}(K(t))$ converges and $\lim_{k\rightarrow\infty}K(\xi_{n_k})=\widehat{K}$,
then $\mathcal{E}(+\infty)
=\lim_{k\rightarrow\infty}\mathcal{E}(K(t_{n_k}))
=\mathcal{E}(\widehat{K})=0$.
Hence, $\lim_{n\rightarrow\infty}\mathcal{E}(K(\xi_n))
=\mathcal{E}(+\infty)=0$.
This is a contradiction.
\qed

\noindent\textbf{Data\ availability\ statements}
Data sharing not applicable to this article as no datasets were generated or analysed during the current study.

\end{document}